\documentclass{amsart}

\usepackage{crimson}
\usepackage{tikz-cd,hyperref,enumerate,tikz-cd}

\theoremstyle{plain}
\newtheorem{thm}{Theorem}[section]

\newtheorem{corollary}{Corollary}[section]

\newtheorem{lemma}{Lemma}[section]
\newtheorem*{msodc}{Motivic Semiorthogonal Decomposition Conjecture (MSODC)}

\theoremstyle{definition}
\newtheorem{defn}{Definition}[section]
\newtheorem{remark}{Remark}[section]

\begin{document}
\title{Note on Motivic Semiorthogonal Decompositions for Elementary Abelian 2-Group Actions}
\author[B. Lim]{Bronson Lim}
\address{BL: Department of Mathematics \\ California State University San
  Bernardino \\ San Bernardino, CA 92407, USA} 
  \email{bronson.lim@csusb.edu}
  
\maketitle

\begin{abstract}
	Let \(\mathcal{X}\) be a smooth Deligne-Mumford stack which is generically a scheme and has quasi-projective coarse moduli. If \(\mathcal{X}\) has elementary Abelian 2-group stabilizers and the coarse moduli of the inertia stack is smooth, we show there exists a semiorthogonal decomposition of the derived category of \(\mathcal{X}\) where the pieces are equivalent to the derived category of the components of the coarse moduli of the inertia stack.
\end{abstract}

\section{Introduction}

\subsection{Motivation}

Let \(X\) be a smooth quasi-projective variety over an algebraically closed field of characteristic zero. Suppose \(G\) is a finite group acting effectively on \(X\). For any additive invariant \(E\) of \(\mathcal{D}[X/G]\) there is a decomposition
\[
	E(\mathcal{D}[X/G]) = \bigoplus\limits_{g\in G/\sim}E(X^g)^{C(g)}
\]
where \(G/\sim\) is the set of conjugacy classes, see \cite[Remark 1.26]{tvdb-orbifold}. Moreover, if each of the quotients \(\bar{X}^g = X^g/C(g)\) are smooth, pullback induces an isomorphism
\[
	E(X^g)^{C(g)}\cong E(\bar{X}^g).
\]
In \cite[Conjecture A]{pvdb-equivariant}, the authors conjecture this additive decomposition admits a categorification at the derived or dg level: 

\begin{msodc}
	Suppose a finite group acts effectively on a quasi-projective variety \(X\) such that each quotient \(\bar{X}^g = X^g/C(g)\) is smooth. Then there is a total order on \(G/\sim\), say \(g_1>g_2>\cdots >g_r\) and a semiorthogonal decomposition
	\[
		\mathcal{D}[X/G] = \langle \mathcal{D}(\bar{X}^{g_1}),\ldots, \mathcal{D}(\bar{X}^{g_r})\rangle.
	\]
\end{msodc}

Although there are many known instances of this conjecture, see \cite[Section 1]{L-R-abelian}, it is worth noting that this conjecture is still open in dimensions larger than 1. Even in the case of surfaces, the \'etale local case \([\mathbb{A}^2/G]\) with \(G\) a rank 2 complex reflection group is open for many groups including \(G(m,p,2)\) with \(p>1\) and each of the exotic rank 2 complex reflection groups.

\subsection{Main Result}

In this note, we give an affirmative answer in the case \(G\) is an elementary Abelian 2-group. For such stacks the inclusions of fixed loci are \(G\)-equivariant and the components of the semiorthogonal decomposition can be canonically ordered by dimension. Further, we show that the conjecture extends to smooth Deligne-Mumford stacks with appropriate conditions.

\begin{thm}
	Let \(\mathcal{X}\) be a smooth Deligne-Mumford stack such that:
	\begin{enumerate}
		\item there is an open dense subscheme \(U\subset\mathcal{X}\);
		\item the coarse moduli \(\pi\colon\mathcal{X}\to X\) is quasi-projective;
		\item the connnected components \(X=X_1,\ldots,X_r\) of the coarse moduli space of the inertia stack \(\mathcal{IX}\) are smooth;
		\item for any \(x\in\mathcal{X}\) the stabilizer group \(\mathrm{st}(x)\) is an elementary Abelian 2-group, i.e. is isomorphic to \(\mu_2^k\) for some \(k\);
	\end{enumerate}
	Then, up to permuting the indices, there is a semiorthogonal decomposition of the form
	\[
		\mathcal{D(X)} = \langle \mathcal{D}(X_1),\ldots,\mathcal{D}(X_r)\rangle,
	\]
	which is \(\mathcal{D}(X)\)-linear.
	\label{thm:main}
\end{thm}

Conditions (1), (2), and (3) are the global analogues of the Zariski local case of the conjecture. 

Previous work towards Theorem \ref{thm:main} can be found in \cite[Example 6.2.4]{potter-thesis} where the case of surfaces and iterated root stacks are considered. Our result is a generalization of his both in terms of dimension and refining the contribution from the canonical stack. 

Moreover, our result is a step towards a complete understanding of the Conjecture in the case of Abelian groups as well as an attempt to globally formulate the Motivic Semiorthogonal Decomposition Conjecture. In particular, we see that when a canonical order in the \'etale local case is provided, it is easier to work in the full generality of stacks rather than with \(G\)-equivariant objects. We expect, in the general Abelian case, to have additional contributions from the normal bundle on each component of the inertia stack or various thickenings of the embedding functors. We leave this for future work.

\subsection{Outline of Proof}
We first work \'etale locally. We then use the \'etale local picture to prove the global case, which is Theorem \ref{thm:main}. We then apply the global case to prove the motivic semiorthogonal decomposition conjecture for elementary Abelian 2-groups.

\subsection{Outline of Paper}

In Section \ref{sec:prelims} we recall preliminary information on semiorthogonal decompositions and fully-faithful functors. In Section \ref{sec:main} we prove the Theorem \ref{thm:main}. We then show it implies the motivic semiorthogonal decomposition conjecture in the case of elementary Abelian 2-groups and give an application to quadrics.

\subsection{Conventions}

We work over an algebraically closed field \(k\) of characteristic zero but the results should hold for tame stacks. For \(\mathcal{X}\) a scheme or stack we denote the bounded derived category of coherent sheaves on \(\mathcal{X}\) by \(\mathcal{D(X)}\). All functors are assumed to be derived unless stated otherwise.

\section{Preliminaries}
\label{sec:prelims}

We recall preliminary material on derived categories of coherent sheaves and semiorthogonal decompositions. The standard reference is \cite{kuznetsov-sod}.

\subsection{Semiorthogonal Decompositions}

\begin{defn}
	Let \(\mathcal{T}\) be a triangulated category. Suppose \(\mathcal{A,B}\) are a pair of triangulated subcategories of \(\mathcal{T}\) such that
	\begin{enumerate}
		\item[(SO)] for each \(a\in \mathcal{A}\) and \(b\in\mathcal{B}\) we have
		\[
			\mathrm{Hom}_{\mathcal{T}}(b,a) = 0;
		\]
		\item[(D)] for each \(t\in \mathcal{T}\), there exists \(a\in\mathcal{A}\) and \(b\in\mathcal{B}\) and an exact triangle
		\[
			b\to t\to a\to b[1].
		\]
	\end{enumerate}
	Then we say \(\mathcal{T} = \langle \mathcal{A,B}\rangle\) is a \textbf{semiorthogonal decompostion} of \(\mathcal{T}\).
\end{defn}

This definition can be iterated to get semiorthogonal decompositions with \(n\) components:
\[
	\mathcal{T} = \langle \mathcal{A}_1,\ldots,\mathcal{A}_n\rangle.
\]

Now suppose \(X\) is a smooth quasi-projective variety and \(G\) is a finite group acting on \(X\). Then the derived category of the stack \([X/G]\) is triangulated. The following definition is motivated by the MSODC.

\begin{defn}
	A \textbf{motivic semiorthogonal decomposition} of \(\mathcal{D}[X/G]\) is the data of
	\begin{itemize}
		\item A total order \(\lambda_1>\lambda_2>\cdots>\lambda_r\) on conjugacy classes;
		\item Fully-faithful embedding functors
		\[
			\Phi_{\lambda_i}\colon \mathcal{D}(X^{\lambda_i}/C_G(\lambda_i))\hookrightarrow \mathcal{D}[X/G]
		\]
		for each \(\lambda_i\in G/\sim\) which are linear over \(\mathcal{D}(X/G)\);
		\item a semiorthogonal decomposition
		\[
			\mathcal{D}[X/G] = \langle C_1,\ldots,C_r\rangle
		\]
		where \(C_i\) is the image of \(\Phi_{\lambda_i}\).
	\end{itemize}
\end{defn}

\subsection{Bondal-Orlov Fully-Faithfulness Criterion}

The original Bondal-Orlov Fully-Faithfulness Criterion was a statement strictly about functors between derived categories of smooth projective varieties. The following is a more general version, see \cite[Theorem 2.7.1]{L-equivariant}.

\begin{thm}[Bondal-Orlov Fully-Faithfulness Criterion]
	Let \(X\) be a smooth quasi-projective variety and \(\Phi\colon\mathcal{D}(X)\to\mathcal{T}\) be a functor to a triangulated category that has a right adjoint. Then \(\Phi\) is fully-faithful provided for each pair of closed points \(x,y\in X\) we have
	\begin{itemize}
		\item If \(x\neq y\), then \(\mathrm{Hom}_{\mathcal{T}}(\Phi(\mathcal{O}_x),\Phi(\mathcal{O}_y)[i]) = 0\) for all \(i\), i.e. \(\Phi\) separates points;
		\item If \(x =y\), then
		\[
			\mathrm{Hom}_{\mathcal{T}}(\Phi(\mathcal{O}_x),\Phi(\mathcal{O}_x)[i]) = \begin{cases} \mathbb{C} & i = 0 \\ 0 & i\notin[0,\dim(X)] \end{cases},
		\]
		i.e. \(\Phi\) induces an isomorphism on tangent spaces.
	\end{itemize}
	Then \(\Phi\) is fully-faithful.
\end{thm}

\section{Proof of Main Result}
\label{sec:main}

\subsection{The \'etale local case.}

Set \(X = \mathbb{A}^n\) and \(G = \mu_2\times\cdots\times\mu_2= \mu_2^k\) where the \(i\)th copy of \(\mu_2\) is generated by \(\lambda_i\) and \(k\leq n\). For a binary multi-index \(I = (i_1,\ldots,i_k)\), we set \(\lambda_I = (\lambda_1^{i_1},\ldots,\lambda_k^{i_k})\). We define and action of \(G\) on \(X\) by
\[
\lambda_1^{i_1}\cdots\lambda_k^{i_k} \cdots (x_1,\ldots,x_k,x_{k+1},\ldots,x_n) = (\lambda_1^{i_1}x_1,\ldots,\lambda_k^{i_k}x_k,x_{k+1},\ldots,x_n).
\]

Define \(|I|\) to be the number of nontrivial entries of \(I\). Define a partial order on conjugacy classes by defining \(I>J\) if \(|I|<|J|\). Then define a total order on conjugacy classes by refining the partial order \(>\) to any total order.

Set \(X^I = X^{\lambda_I}\) to be the fixed locus of the action of \(\lambda_I\) on \(X\). We have an isomorphism \(X^I\cong \mathbb{A}^{n_I}\) where \(n_I = n - |I|\) and hence an isomorphism \(\bar{X}^I\cong\mathbb{A}^{n_I}\) as well. Let \(\pi_I\colon X^I\to\bar{X}^I\) be the quotient mapping and \(\iota_I\colon X^I\hookrightarrow \mathbb{A}^n\) be the \(\mu_2^k\)-equivariant inclusion.

Now set
\[
	\Phi_I= \iota_{I\ast}\circ\pi_I^\ast\colon\mathcal{D}(\bar{X}^I)\to \mathcal{D}[\mathbb{A}^n/\mu_2^k].
\]
Note that \(\Phi_i\) has both a left and right adjoint using the standard (equivariant) adjunction formula.

\begin{thm}
	The total order on binary multi-indices induces a motivic semiorthogonal decomposition
	\[
		\mathcal{D}[\mathbb{A}^n/G] = \langle \Phi_I\mathcal{D}(\bar{X}^I)\rangle_{I}
	\]
	where for each \(I>J\), the subcategory \(\langle \mathcal{D}(\bar{X}^I),\mathcal{D}(\bar{X}^J)\rangle\) is semiorthogonal. That is, the subcategories are orded by decreasing dimension.
	\label{thm:etale-local}
\end{thm}

\begin{proof}
	This is essentially \cite[Example 4.3.3]{L-P-divisor}. We recall it now for completeness. We have an equivalence
	\[
		[\mathbb{A}^n/\mu_2^k]\cong [\mathbb{A}^1/\mu_2]^k\times \mathbb{A}^{n-k}.
	\]
	Let \(\pi\colon \mathbb{A}^1\to\mathbb{A}^1\) be the quotient mapping given by \(x\mapsto x^2\). Let \(\iota\colon\{0\}\hookrightarrow \mathbb{A}^1\) denote the \(\mu_2\)-equivariant inclusion of the origin. Then the derived category of the stack \([\mathbb{A}^1/\mu_2]\) admits a semiorthogonal decomposition of the form:
	\[
		\mathcal{D}[\mathbb{A}^1/\mu_2] = \langle \pi^\ast\mathcal{D}(\mathbb{A}^1),\iota_\ast\mathcal{D}(0)\rangle.
	\]
	For each \(i = 1,\ldots, k\) set \(\pi_i\colon \mathbb{A}^1_{x_i}\to \mathbb{A}^1_{x_i^2}\) to be the quotient mapping. Then
	\[
		(\pi_{i_1},\ldots,\pi_{i_r})\colon \mathbb{A}^r_{x_I} := \mathbb{A}^1_{x_{i_1}}\times\cdots \mathbb{A}^1_{x_{i_r}}\to \mathbb{A}^r_{x_I^2}:=\mathbb{A}^1_{x_{i_1}^2}\times\cdots \mathbb{A}^1_{x_{i_r}^2}
	\]
	is the corresponding quotient mapping. Let \(I\) be the binary multi-index corresponding to \(i_1,\ldots,i_r\) so that \(\pi_I = (\pi_{i_1},\ldots,\pi_{i_r})\). Passing to a dg lift, we have a quasi-equivalence
	\begin{align*}
		\mathcal{D}[\mathbb{A}^n/\mu_2^k] &\cong \mathcal{D}[\mathbb{A}^1/\mu_2]^{\otimes k}\otimes\mathcal{D}(\mathbb{A}^{n-k}) \\
		&\cong \left(\bigotimes_{i=1}^k\langle \pi_i^\ast\mathcal{D}(\mathbb{A}_{x_i}^1),\iota_\ast(0)\rangle\right)\otimes \mathcal{D}(\mathbb{A}^{n-k}) \\
		&\cong \left\langle \mathcal{D}(\mathbb{A}^{n-k}/\mu_2^k), \bigoplus_{|I|=1}\Phi_I(\mathcal{D}(\mathbb{A}^{n-k-1}_{x_I^2})),\ldots,\mathcal{D}(0)\right\rangle \otimes \mathcal{D}(\mathbb{A}^{n-k})\\
		&\cong \left\langle \mathcal{D}(\mathbb{A}^n/\mu_2^k), \bigoplus_{|I|=1}\Phi_I(\mathcal{D}(\bar{X}^I)),\ldots,\mathcal{D}(X^{(1,\ldots,1)})\right\rangle.
	\end{align*}
	This completes the proof.
\end{proof}

\begin{corollary}
	Let \(U\subset\mathbb{A}^n\) be a \(\mu_2^k\)-invariant open subset. Then restriction of \(\Phi_I\) defines a motivic semiorthogonal decomposition of \(\mathcal{D}[U/\mu_2^k]\).
	\label{cor:msod-open}
\end{corollary}

\begin{proof}
	By \cite[Corollary 2.3.8]{L-P-divisor}, the motivic semiorthogonal decompositoin for \([\mathbb{A}^n/\mu_2^k]\) restricts to \([U/\mu_2^k]\) for any \(\mu_2^k\)-invariant open subset \(U\).
\end{proof}

\subsection{Global Case}

Let \(\mathcal{X}\) be a smooth DM stack with smooth and quasi-projective coarse moduli. Additionally, suppose \(\mathcal{X}\) has elementary Abelian 2-group stabilizers and is generically a scheme. We consider the following diagram
\[
\begin{tikzcd}
{} & \mathcal{IX}\cong \coprod_i\mathcal{X}_i \ar{d}{\pi} \ar{dr}{\rho} & {} \\
\bar{\mathcal{X}}_i\ar{r}{\iota_i} & \coprod_i \bar{\mathcal{X}}_i & \mathcal{X}
\end{tikzcd}
\]
where \(\mathcal{X}_i\) are the connected components of \(\mathcal{IX}\) and \(\mathcal{X}_i\to\bar{\mathcal{X}}_i\) their coarse moduli which we assume to be smooth. We define Fourier-Mukai functors
\[
	\Phi_i = \rho_\ast\circ\pi^\ast\circ\iota_{i\ast}\colon \mathcal{D}(\bar{\mathcal{X}}_i)\to \mathcal{D(X)}.
\]
We note that \(\Phi_i\) has both a left and right adjoint using the standard formulas.

Define a dimensional partial order on the set of connected components, i.e. \(i>j\) if \(\dim(\bar{\mathcal{X}}_i)>\dim(\bar{\mathcal{X}}_j)\). Take any refinement to a total order and relabel so that the total order is the same as the natural order \(1>2>\cdots>r-1>r\). We will show \(\Phi_i\) define a semiorthogonal decomposition of \(\mathcal{D}(\mathcal{X})\).

First notice that if \(x\in\bar{\mathcal{X}}_i\) is a closed point, then we can consider its image \(\bar{x}\in \bar{\mathcal{X}}\) and the unique lift \(\tilde{x}\in\mathcal{X}\). There is an inertia-preserving \'etale mapping \(\rho\colon [U/\mathrm{st}(\tilde{x})]\to\mathcal{X}\), see \cite[Theorem 2.12]{olsson-homstacks}, where \(U\subseteq\mathbb{A}^n\) is open. In particular, the diagram
\[
	\begin{tikzcd}
	\mathcal{I}{[}U/\mathrm{st}(\tilde{x}){]}\ar{r}{\rho'} \ar{d}{\sigma'} & {[}U/\mathrm{st}(\tilde{x}){]} \ar{d}{\sigma} \\
	\mathcal{IX} \ar{r}{\rho} & \mathcal{X}
	\end{tikzcd}
\]
is Cartesian and the vertical arrows are \'etale and surjective in a neighborhood of \(\tilde{x}\), i.e. those connected components of \(\mathcal{IX}\) that map to \(\bar{x}\in\bar{\mathcal{X}}\). Without loss of generality, we will assume \(\sigma'\) is surjective in a neighborhood of \(\tilde{x}\) as our computations are local. Lastly, the action of \(\mathrm{st}(\tilde{x})\) on \(U\) is effective since \(\mathcal{X}\) is generically a scheme.

Consider now the extended commutative diagram:
\begin{equation}
	\begin{tikzcd}
	U^\mu/\mathrm{st}(\tilde{x})\ar{r}{\iota_\mu}\ar{d} & \coprod_{\lambda\in\mathrm{st}(\tilde{x})} U^\lambda/\mathrm{st}(\tilde{x})\ar{d}{\bar{\sigma}'} & \mathcal{I}{[}U/\mathrm{st}(\tilde{x}){]}\ar{r}{\rho'} \ar{d}{\sigma'} \ar{l}{\pi'} & {[}U/\mathrm{st}(\tilde{x}){]} \ar{d}{\sigma} \\
	\bar{\mathcal{X}}_j \ar{r}{\iota_j} & \coprod_i\bar{\mathcal{X}}_i & \mathcal{IX} \ar{r}{\rho}\ar{l}{\pi} & \mathcal{X}
	\end{tikzcd}
\end{equation}
where all vertical arrows are \'etale and the leftmost map is just restriction of \(\bar{\sigma}'\). The functors
\[
	\Phi_\mu = \rho'_\ast\circ (\pi')^\ast\circ \iota_{\mu\ast}\colon \mathcal{D}(U^\mu/\mathrm{st}(\tilde{x}))\to \mathcal{D}[U/\mathrm{st}(\tilde{x})]
\]
are fully-faithful for each \(\mu\in\mathrm{st}(\tilde{x})\) and define a motivic semiorthogonal decomposition ordered by decreasing dimension by Corollary \ref{cor:msod-open}.

\begin{lemma}
	For each \(i\) the functor \(\Phi_i\colon\mathcal{D}(\bar{\mathcal{X}}_i)\to\mathcal{D}(\mathcal{X})\) is fully-faithful.
	\label{lem:ff}
\end{lemma}

\begin{proof}
	We apply the Bondal-Orlov fully-faithfulness criterion. Since \(\Phi_i\) separates points it suffices to check that
	\[
		\mathrm{Hom}_\mathcal{X}(\Phi_i(\mathcal{O}_x),\Phi_i(\mathcal{O}_x)[0])\cong \mathbb{C}
	\]
	and vanishes outside \([0,\dim(\bar{\mathcal{X}}_i)]\).
	
	For both claims consider the commutative diagram (1). Pick a lift \(z\in(\bar{\sigma}')^{-1}(x)\), then we have \(\bar{\sigma}'_\ast\mathcal{O}_z\cong \mathcal{O}_x\). Since the vertical maps are \'etale the base change map is an isomorphism for the middle square:
	\[
		\pi^\ast\circ \bar{\sigma}'_\ast\cong \sigma'_\ast\circ(\pi')^\ast.
	\]
	Thus we have
	\begin{align*}
		\Phi_i(\mathcal{O}_x) &\cong \Phi_i(\bar{\sigma'}_\ast\mathcal{O}_z) \\
		&\cong (\rho_\ast\circ\pi^\ast\circ\iota_{j\ast})(\bar{\sigma'}_\ast\mathcal{O}_z) \\
		&\cong (\rho_\ast\circ\pi^\ast)((\bar{\sigma'}_\ast\circ\iota_{\mu\ast})(\mathcal{O}_z)) \\
		&\cong (\rho_\ast\circ(\pi^\ast\circ\bar{\sigma'}_\ast)\circ\iota_{\mu\ast})(\mathcal{O}_z)) \\
		&\cong (\rho_\ast\circ(\sigma'_\ast\circ(\pi')^\ast)\circ\iota_{\mu\ast})(\mathcal{O}_z)) \\
		&\cong (\sigma_\ast\circ(\rho')^\ast)\circ((\pi')^\ast\circ\iota_{\mu\ast})(\mathcal{O}_z)) \\
		&\cong \sigma_\ast(\Phi_\mu(\mathcal{O}_z)) \\
	\end{align*}
	Now since \(\sigma\) is \'etale, by Grothendieck duality we have:
	\begin{align*}
		\mathcal{H}om_{\mathcal{X}}(\Phi_i(\mathcal{O}_x),\Phi_i(\mathcal{O}_x)) &\cong \mathcal{H}om_{\mathcal{X}}(\sigma_\ast(\Phi_\mu(\mathcal{O}_z)),\sigma_\ast(\Phi_\mu(\mathcal{O}_z))) \\
		&\cong \sigma_\ast\mathcal{H}om_{[U/\mathrm{st}(\tilde{x})]}(\Phi_\mu(\mathcal{O}_z),\sigma^\ast\sigma_\ast(\Phi_\mu(\mathcal{O}_z))) \\
		&\cong \sigma_\ast\mathcal{H}om_{[U/\mathrm{st}(\tilde{x})]}\left(\Phi_\mu(\mathcal{O}_z),\bigoplus_{y\in\sigma^{-1}(\tilde{x}=\sigma(y))}(\Phi_\mu(\mathcal{O}_y))\right) \\
		&\cong \sigma_\ast\mathcal{H}om_{[U/\mathrm{st}(\tilde{x})]}(\Phi_\mu(\mathcal{O}_z),\Phi_\mu(\mathcal{O}_z)) \\
	\end{align*}
	where the last line follows since the first argument is supported at \(z\). This implies
	\[
		\mathrm{Hom}_{\mathcal{X}}(\Phi_i(\mathcal{O}_x),\Phi_i(\mathcal{O}_x)) \cong \mathrm{Hom}_{[U/\mathrm{st}(\tilde{x})]}(\Phi_\mu(\mathcal{O}_z),\Phi_\mu(\mathcal{O}_z))
	\]
	and the claims follow.
\end{proof}

\begin{remark}
If we do not assume \(\bar{\mathcal{X}}\) is smooth then the proof of Lemma \ref{lem:ff} carries over verbatim provided \(\bar{\mathcal{X}}_i\) is smooth and the image of \(\bar{\mathcal{X}}_i\to\bar{\mathcal{X}}\) avoids \(\mathrm{Sing}(\bar{\mathcal{X}})\). We will not need this generality; however, it would be interesting to determine the semi-orthogonal complement in this case. In the case of a surface, it is the derived category of the minimal resolution of the coarse moduli by the Bridgeland-King-Reid theorem.
\end{remark}

\begin{lemma}
	If \(i>j\), then \(\langle \Phi_i\mathcal{D}(\bar{X}_i),\Phi_j\mathcal{D}(\bar{X}_j)\rangle\) is semiorthogonal.
	\label{lem:so}
\end{lemma}

\begin{proof}
	Since our total order is given by dimension, which is the same as the local description, semiorthogonality follows from a similar computation to the proof of Lemma \ref{lem:ff}.
\end{proof}

We can now prove the main result.

\begin{proof}[Proof of Theorem \ref{thm:main}]
	From Lemmas \ref{lem:ff} and \ref{lem:so}, it suffices to prove generation. Again, we can check this \'etale locally and the claim follows.
\end{proof}

\subsection{A Zariski Local Example}

Before we prove the MSODC for elementary Abelian 2-groups, we show that the ordering we are given from Theorem \ref{thm:main} may not match the ordering of the MSODC. We may have to, as this example shows, mutate the components to match.

Consider the action of \(G= \mu_2^2\) on \(X=\mathbb{P}^2\) by multiplication of the first two homogeneous coordinates:
\[
	( (-1)^a,(-1)^b)\cdot [x:y:z] = [(-1)^ax:(-1)^by:z].
\]

The nontrivial fixed loci are
\begin{itemize}
	\item \(X^{(-1,1)} = \{p = [1:0:0]\}\cup V(x)\);
	\item \(X^{(1,-1)} = \{q = [0:1:0]\}\cup V(y)\);
	\item \(X^{(-1,-1)} = \{r = [0:0:1]\}\cup V(z)\).
\end{itemize}
The coarse moduli of each fixed locus is isomorphic to the disjoint union of a point and a line. The components of the coarse moduli of the inertia stack are:
\begin{itemize}
	\item Three points: \(p,q,r\).
	\item Three lines: \(\mathbb{P}^1_x, \mathbb{P}^1_y, \mathbb{P}^1_z\).
	\item One plane: \(\mathbb{P}^2\).
\end{itemize}
By Theorem \ref{thm:main} we have a semiorthogonal decomposition of the form
\[
	\mathcal{D}[\mathbb{P}^2/\mu_2^2] = \langle \mathcal{D}(\mathbb{P}^2), \mathcal{D}(\mathbb{P}^1_x), \mathcal{D}(\mathbb{P}^1_y), \mathcal{D}(\mathbb{P}^1_z), \mathcal{D}(p), \mathcal{D}(q),\mathcal{D}(r)\rangle.
\]
Since the images of the derived categories of the lines (resp. points) are mutually orthogonal, we can interchange their order. However, the derived categories of points that appear are not mutually orthogonal to the derived categories of lines. Nevertheless, we can mutate \(\mathcal{D}(p)\) through \(\langle \mathcal{D}(\mathbb{P}^1_y),\mathcal{D}(\mathbb{P}^1_z)\rangle\) and then mutate \(\mathcal{D}(q)\) through \(\mathcal{D}(\mathbb{P}^1_z)\) to  get
\[
	\mathcal{D}[\mathbb{P}^2/\mu_2^2] = \langle \mathcal{D}(\mathbb{P}^2), \mathcal{D}(\mathbb{P}^1_x), \mathcal{D}(p), \mathcal{D}(\mathbb{P}^1_y), \mathcal{D}(q), \mathcal{D}(\mathbb{P}^1_z), \mathcal{D}(r)\rangle.
\]

\subsection{Application to the MSODC}

As our primary application of Theorem \ref{thm:main}, we prove the motivic semiorthogonal decomposition conjecture for elementary Abelian 2-groups. 

\begin{corollary}
	The motivic semiorthogonal decomposition conjecture holds for elementary Abelian 2-groups.
	\label{cor:msod}
\end{corollary}

\begin{proof}
	Let \(G = \mu_2^n\) act effectively on a smooth quasi-projective variety \(X\) such that each of the coarse quotients \(X^\lambda/G\) are smooth for each \(\lambda\in G\). Thus by Theorem \ref{thm:main}, there is a semiorthogonal decomposition of \(\mathcal{D}[X/G]\) with components equivalent to the derived categories of the components of the coarse moduli of the inertia stack. These components are precisely the components of \(X^\lambda/G\). Since mutating the provided \(\mathcal{D}(X/G)\)-linear semiorthogonal decomposition stays \(\mathcal{D}(X/G)\)-linear, we conclude \(\mathcal{D}[X/G]\) has a motivic semiorthogonal decomposition.
\end{proof}

\subsection{Application to Fermat Quadrics}

Consider the normal form to a smooth quadric, i.e. a smooth Fermat quadric
\[
	f(x_1,\ldots,x_{n+1}) = \sum_{i=1}^{n+1} x_i^2.
\]
Then the vanishing locus of \(f\) defines a smooth quadric \(Q\subset \mathbb{P}^n\). The group \(\mu_2^n\) acts on \(\mathbb{P}^n\) by multiplication on the first \(n\) entries in the natural way. This action leaves \(Q\)-invariant.  

\begin{thm}
	Let \(\mu_2^n\) act on \(Q\) as above. Then the derived category \(\mathcal{D}[Q/\mu_2]\) has a full exceptional collection.
	\label{thm:quadric}
\end{thm}

\begin{proof}
For each binary multi-index \(I\), again let \(|I|\) again denote the number of non-trivial entries. If \(|I| = 1\), then the corresponding fixed locus \(Q^I\) can be identified with the restriction of \(Q\) to the remaining \(n\) variables, i.e. a smooth Fermat quadric in \(\mathbb{P}^{n-1}\). If \(1<|I|<n-1\), \(Q^I\) is the union of two Fermat quadrics given by restriction to the non-trivial variables and to the trivial variables. If \(|I|=n-1\), then \(Q^I\) is the restriction of \(Q\) to the first \(n\) variables which is again a smooth Fermat quadric in \(\mathbb{P}^{n-1}\). Finally, each of the quotients \(Q^I/\mu_2^n\) is isomorphic to either \(\mathbb{P}^k\) or a disjoint union of two projective spaces. By Corollary \ref{cor:msod}, \(\mathcal{D}[Q/\mu_2^n]\) has a motivic semiorthogonal decomposition where each piece is a copy of the derived category of a projective space. Thus \(\mathcal{D}[Q/\mu_2^n]\) has a full exceptional collection.
\end{proof}

\begin{remark}
	It's clear how to generalize the Zariski local example of \([\mathbb{P}^2/\mu_2^2]\) above to see that the derived category of \([\mathbb{P}^n/\mu_2^n]\) possesses a motivic semiorthogonal decomposition. Using \cite[Theorem 1.2.1]{L-P-divisor} we can pass this motivic semiorthogonal decomposition to \(Q\). The resulting motivic semiorthogonal decomposition is the same as the one of Theorem \ref{thm:quadric}.
\end{remark}

\bibliographystyle{amsalpha}
\bibliography{msod-abelian}

\end{document}